\documentclass[11pt]{article}
\usepackage{amsfonts}
\usepackage{lscape}
\usepackage{latexsym,amsmath}
\usepackage{amssymb,array}
\makeatletter \oddsidemargin 0in \evensidemargin 0in \textwidth
16cm \RequirePackage[dvips]{graphicx} \textheight 20cm
\newtheorem{t1}{Theorem}[section]
\newtheorem{p1}{Proposition}[section]

\newtheorem{r1}{Remark}[section]

\newtheorem{ex}{Example}[section]

\begin{document}
\title{On Weighted Generalized Entropy for Double Truncated Distribution}
\author{Shivangi Singh\\\and Chanchal Kundu\footnote{Corresponding
author e-mail: chanchal$_{-}$kundu@yahoo.com/ckundu@rgipt.ac.in}
\and\and
Department of Mathematics\\
Rajiv Gandhi Institute of Petroleum Technology\\
Jais 229 304, U.P., India}
\date{March, 2020}
\maketitle
\begin{abstract}
 The notion of weighted Renyi's entropy for truncated random variables has recently been proposed in the information-theoretic literature. In this paper, we introduce a generalized measure of it for double truncated distribution, namely weighted generalized interval entropy (WGIE), and study it in the context of reliability analysis. Several properties, including monotonicity, bounds and uniqueness of WGIE are investigated. Moreover, a simulation study is carried out to demonstrate the performance of the estimates of the proposed measure using simulated and real data sets. The role of WGIE in reliability modeling has also been investigated for a real-life problem.
\end{abstract}
\section{Introduction and some preliminary results}
In the literature, the notion of weighted distribution was introduced by Fisher (1934) and later explored by Rao (1965) in connection with modeling statistical data where the usual practice of employing standard distributions for the purpose was not found appropriate. For instance, in many real life situations when an investigator collects a sample of observations of any practical event, the standard distributions may not be fitted due to various reasons such as non-observe ability of some events or damage caused to the original observations. These sampling situations can be modeled by using weighted distributions. Weighted distributions arise when the observations generated from a stochastic process are recorded with some weight function and are frequently studied in areas; such as survival analysis, reliability, analysis of family data, bio-medicine, forestry, ecology and survey sampling, to mention a few.\\
\hspace*{.2in} Shannon (1948) entropy is a very important and well-known concept in the field of information theory, statistics, data compression, engineering sciences, especially in communication engineering. It is a shift independent measure and gives equal importance or weight to the occurrence of every event. However, in some practical situations, such as reliability or neurobiology, a shift-dependent measure of uncertainty is desirable. To this aim, Belis and Guia\c{s}u (1968) and later Di Crescenzo and Longobardi (2006) considered the notion of weighted entropy. Let $X$ be an absolutely continuous nonnegative random variable with distribution function $F$ and probability density function $f$. Then the weighted entropy is defined as
\begin{eqnarray}\label{eq1.1}
H^{w}(X)=-\int_0^\infty x f(x)\log f(x)dx=-\int_0^\infty\int_y^\infty f(x)\log f(x)dydx.
\end{eqnarray}
As pointed out by Belis and Guia\c{s}u (1968) that the occurrence of an event removes a double uncertainty: the quantitative one, related to the probability with which it occurs, and the qualitative one, related to its utility for the attainment of the goal or to its significance with respect to a given qualitative characteristic. The factor $x$, in the integral on the right-hand-side of (\ref{eq1.1}), may be viewed as a weight linearly emphasizing the occurrence of the event $\{X=x\}$. This yields a length biased shift-dependent information measure assigning greater importance to larger values of $X$. The use of weighted entropy (\ref{eq1.1}) is also motivated by the need, arising in various communication and transmission problems, of expressing the usefulness of events with an information measure. An important feature of the human visual system is that it can recognize objects in a scale and translation invariant manner. Achieving this desirable behavior using biologically realistic networks is a challenge (cf. Wallis, 1996). Indeed, knowing that a device fails to operate, or a neuron fails to release spikes in a given time-interval, yields a relevantly different information from the case when such an event occurs in a different equally wide interval. In some cases we are thus led to resort to a shift-dependent information measure that, for instance, assigns different measures to such distributions.\\
\hspace*{.2in} Recently, based on the idea of weighted entropy, Das (2017) introduced the concept of first and second kind weighted entropies of order $\alpha$ and studied their properties in the context of left/right truncated random variable. In the same vein of second kind weighted entropy of order $\alpha$, Nourbakhsh and Yari (2017) introduce weighted Renyi's (1961) entropy
\begin{equation}\label{eq1.2}
 H_\alpha^{w}(X)=\frac{1}{1-\alpha}\log\int_0^\infty \left(x f(x)\right)^{\alpha}dx;~~for~0<\alpha\neq1,
\end{equation}
where the integral is finite. The factor $x$ in the right hand side integral yields a shift-dependent information measure assigning greater importance to larger values of the random variable $X$. In recent years, the study of Renyi's entropy based on the weighted notion has attracted the attention of a number of researchers, such as Sekeh et al. (2014), Das (2017) and Rajesh et al. (2017).\\
\hspace*{.2in} Recall that one important generalization of the Renyi entropy is the Varma (1966) entropy. Based on Varma's entropy, a two-parametric generalization of Shannon entropy, called generalized entropy of order $(\alpha,\beta)$, is given by
\begin{equation}\label{eq1.4}
H_{\alpha,\beta}(X)=\frac{1}{\beta-\alpha}\log\int_0^\infty f^{\alpha+\beta-1}(x)dx,~~for~\beta-1<\alpha<\beta,~\beta\geq 1.
\end{equation}
It plays a vital role as a measure of complexity and uncertainty in different areas such as physics, electronics and engineering to describe many chaotic systems. For some recent work based on (\ref{eq1.4}) one may refer to Baig and Dar (2008), Kumar and Taneja (2011), Rajesh et al. (2014), Kayal (2015), Kundu (2015) and Minimol (2017), to mention a few. Very recently, Kundu and Singh (2020) discuss the usefulness of (\ref{eq1.4}) over the several other generalized entropy measures introduced in the literature including the role of the parameters $\alpha,\beta$ and how to get their values in practice, the advantage of adding the extra parameter $\beta$. Motivated with the utility of generalized entropy and weighted Renyi's entropy, the concept of weighted generalized entropy has been introduced here. In view of (\ref{eq1.4}) we propose a two-parametric generalization of weighted Renyi's entropy, called weighted generalized entropy, which is defined as
\begin{equation}\label{eq1.3}
H_{\alpha,\beta}^{w}(X)=\frac{1}{\beta-\alpha}\log\int_0^\infty \left(x f(x)\right)^{\alpha+\beta-1}dx,~~for~\beta-1<\alpha<\beta,~\beta\geq 1.
\end{equation}
This new measure is shift-dependent and a generalization of recent weighted Renyi's entropy measure. Intuitively, (\ref{eq1.3}) is interpreted as a measure of uncertainty supplied by a probabilistic experiment depending both on the probabilities of events and on qualitative weights of the possible events. The following example illustrates the importance of qualitative characteristic of information as reflected in the definition of weighted generalized entropy.
\begin{ex}
Let $X$ and $Y$ denote random lifetimes of two components with probability density functions
\begin{equation*}f_X(t)=\left\{\begin{array}{ll}\frac{1+t}{4},
~if~0\leq t<2,\\
0,
\quad otherwise
\end{array}\right.
~~and~~~
f_Y(t)=\left\{\begin{array}{ll}1-\frac{1+t}{4},
~if~0\leq t<2,\\
0,
\quad otherwise.
\end{array}\right.
\end{equation*}
For $\alpha = 0.5$ and $\beta=1.2$ we obtain $H_{\alpha,\beta}(X)= H_{\alpha,\beta}(Y)=0.283991$. But, $H_{\alpha,\beta}^{w}(X)=0.346064$ and $H_{\alpha,\beta}^{w}(Y)=0.0809797$. Hence, even though $H_{\alpha,\beta}(X)= H_{\alpha,\beta}(Y)$, the weighted generalized entropy about the predictability of $X$ by the density function $f_X(t)$ is greater than the predictability of $Y$ by the density function $f_Y(t)$. Nevertheless, the generalized entropies measured from a quantitative point of view, neglecting the qualitative side, fails to make any distinction whatsoever between them. To distinguish them, we must take into account the qualitative characteristic as given in (\ref{eq1.3}).
\end{ex}
\hspace*{.2in} The main objective of our present study is to extend the concept of weighted Renyi's entropy for truncated random variables to weighted generalized interval entropy. The remainder of the paper is arranged as follows. In Section 2 we introduce the notion of weighted generalized interval entropy (WGIE) for doubly truncated random variable and study its properties. We obtain upper and lower bounds for the proposed concept and also discuss its monotonicity. It is shown that the WGIE determines the distribution uniquely. In Section 3, Monte-Carlo simulation is carried out to analyze the behavior of the estimates of WGIE which are validated using simulated and real data sets. Finally, Section 4 concludes the present study with an application of the proposed measure in reliability modeling.
\section{Weighted generalized interval entropy and its properties}
Recently, there has been growing interest to study (weighted) entropy measure for doubly truncated random variable which has far-reaching applications in many areas. Doubly truncated failure time arises if the event time of individual which falls in a specific time interval are only observed. Moreover, in many survival studies for modeling real-life data, information about lifetime between two points is only available. With this motivation, we introduce the notion of weighted generalized interval entropy.\\
\hspace*{.2in} Let us consider a nonnegative absolutely continuous doubly truncated random variable $(X|t_1\leq X\leq t_2)$ where $(t_1,t_2)\in D=\{(u,v)\in\Re_+^2 :F(u)<F(v)\}$. Then the weighted generalized entropy of order $(\alpha,\beta)$ for $X$ at interval $(t_1,t_2)$, termed as weighted generalized interval entropy (WGIE), is given by
\begin{eqnarray}\label{eq2}
H_{\alpha,\beta}^{w}(X; t_1,t_2)&=&\frac{1}{\beta-\alpha}\log\int_{t_1}^{t_2}\left(\frac{xf(x)}{ F(t_2)-F(t_1)}\right)^{\alpha+\beta-1}dx,
\end{eqnarray}
where $\beta-1<\alpha<\beta,~\beta\geq1$. When the system has the age $t_1$, for different values of $(\alpha,\beta)$, $H_{\alpha,\beta} ^{w}(X; t_1,t_2)$ provides the quantitative-qualitative information spectrum of the remaining life of the system until age $t_2$. Clearly, $H_{\alpha,\beta}^{w}(X; 0,\infty)= H_{\alpha,\beta}^{w}(X)$ as given in (\ref{eq1.3}). For $\beta=1$, we get weighted Renyi's interval entropy as studied by Singh and Kundu (2019). Also, for $\beta=1$, $t_2\rightarrow\infty$ and $t_1\rightarrow 0$ we get weighted residual/past Renyi's entropies which are given, respectively, by
 \begin{equation}\label{eq1.7}
 H_\alpha^{w}(X, t)=\frac{1}{1-\alpha}\log\int_t^\infty\left(x \frac{f(x)}{\overline F(t)}\right)^{\alpha}dx
\end{equation}
\begin{equation}\label{eq1.8}
{\rm~and,~}\overline H_\alpha^{w}(X, t)=\frac{1}{1-\alpha}\log\int_0^t\left(x \frac{f(x)}{ F(t)}\right)^{\alpha}dx,
\end{equation}
for $0<\alpha\neq1$ and studied by Nourbakhsh and Yari (2017). The following example clarifies the effectiveness of the weighted generalized interval entropy.
\begin{ex}
Let $X$ and $Y$ denote random lifetimes of two components with probability density functions $f(x)=\frac{x}{2},~x\in(0,2)$ and $g(x)=2(1-x),~x\in(0,1)$, respectively. Since $X$ and $Y$ belong to different domains, the use of weighted generalized entropy (\ref{eq1.3}) to compare them informatively
is not interpretable. The WGIE in interval $(0.5,0.8)$ are $H_{\alpha,\beta}^{w}(X;0.5,0.8)=1.78963$ and $H_{\alpha,\beta} ^{w}(Y;0.5,0.8)=1.34467$, respectively for $\alpha=1.5$ and $\beta=2$. Hence, the weighted generalized interval entropy for $X$ is greater than $Y$ in the interval $(0.5,0.8)$.
\end{ex}
\hspace*{.2in} In Table 1, we compute $H_{\alpha,\beta}^{w}(X; t_1,t_2)$ for some well-known distributions where $\gamma(s,t)=\int_0^{t}x^{s-1}e^{-x}dx$ is lower incomplete gamma function so that $\int_u^{v}t^{s-1}e^{-t}dt=\gamma(s,v)-\gamma(s,u)$.\\
\hspace*{.2in} Now we investigate different properties including monotonicity and bounds of WGIE. Before stating the results, recall that the general failure rate (GFR) functions of a doubly truncated random variable $(X|t_1<X<t_2)$ are defined as $h_i(t_1,t_2)=f(t_i)/\left(F(t_2)-F(t_1)\right),~i=1,2$. First we discuss the monotonicity of $H_{\alpha,\beta}^{w}(X; t_1,t_2)$ in view of generalized Pareto distribution (GPD) which plays an important role in reliability, extreme value theory and other branches of statistics.
\begin{ex}\label{ex3}
Let $X$ follow GPD having survival function $\overline F(x)=(1+\theta x)^{-\frac{1}{\theta}}$ for $x,\theta>0$, or $0<x<-1/\theta,~\theta<0$. For $\theta\rightarrow0$, the GPD becomes standard exponential distribution. When $\theta>0$, the family of GPD reduces to Pareto Type-II distribution or Lomax distribution. Also, for $\theta<0$, it becomes the power distribution. For GPD
\begin{equation*}
H_{\alpha,\beta}^{w}(X; t_1,t_2)=\frac{1}{\beta-\alpha}\log\int_{t_1}^{t_2}\left(\frac{x(1+\theta x)^{-\frac{1}{\theta}-1}}{(1+\theta t_1)^{-\frac{1}{\theta}}-(1+\theta t_2)^{-\frac{1}{\theta}}}\right)^{\alpha+\beta-1}dx,
\end{equation*}
which is increasing in $t_1$ and $t_2$ (keeping the other fixed) for $\theta=0.8$ and $\alpha+\beta<(>)2$, as shown in Figure \ref{fig3}. Here, $t_1=-\log(u)$ and $t_2=-\log(v)$ have been used while plotting the curves so that $H_{\alpha,\beta}^{w}(X; t_1,t_2)=H_{\alpha,\beta}^{w}(u,v)$, say.
\begin{figure}[ht]
\centering
\begin{minipage}[b]{0.4\linewidth}
\includegraphics[height=4.5cm]{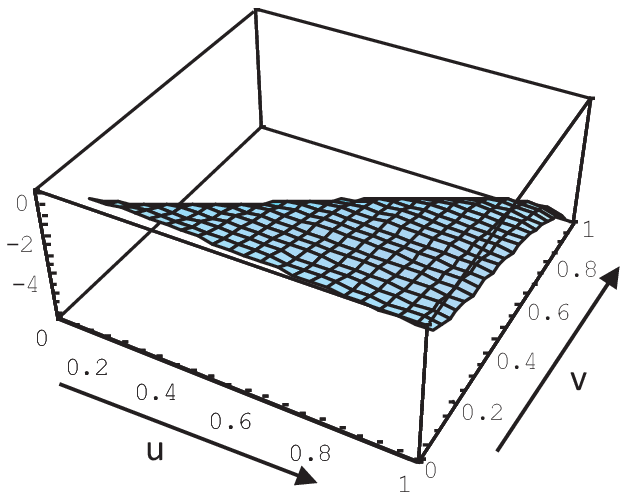}
\centering{$(i)$~Plot of $H_{\alpha,\beta}^{w}(u,v)$ for $\alpha=0.5$ and $\beta=1.2$ against $(u, v)\in(0,1)$}
\end{minipage}
\quad
\begin{minipage}[b]{0.4\linewidth}
\includegraphics[height=4.5cm]{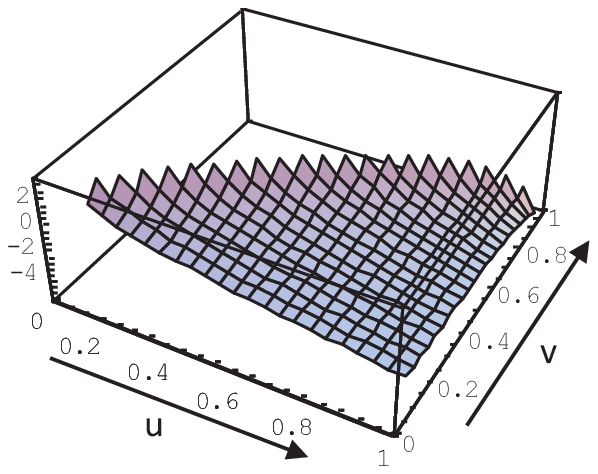}
\centering{$(ii)$~Plot of $H_{\alpha,\beta} ^{w}(u,v)$ for $\alpha=1.2$ and $\beta=2$ against $(u, v)\in(0,1)$}
\end{minipage}
\caption{\label{fig3}Graphical representation of $H_{\alpha,\beta}^{w}(u,v)$ for $\alpha+\beta<(>)2$ (Example \ref{ex3})}
\end{figure}
\end{ex}
\hspace*{.2in} Since the entropy as a measure of uncertainty is expected to decrease when the object's outcome is captured in an interval which is contracting. Recently, Shangri and Chen (2012) give necessary and sufficient condition for the Renyi's interval entropy of an absolutely continuous random variable be an increasing function of interval. An analogous result in the context of weighted Renyi's interval entropy is presented in Singh and Kundu (2018). Following the same, below we have shown that this intuitive monotonicity is also preserved for $H_{\alpha,\beta}^{w}(X; t_1,t_2)$.
\begin{t1}\label{th3.11}
Let $X$ be an absolutely continuous random variable with density function $f(x)$ and twice-differentiable cumulative distribution function $F(x)$. If $F(x)$ is log-concave then, for $\alpha+\beta<2$, $H_{\alpha,\beta}^{w}(X; t_1,t_2)$ is\\
(i) increasing in $t_2$ if $F(t_1)=0$;\\
(ii) a partially increasing function of interval $(t_1,t_2)$.
\end{t1}
Proof: (i) For $F(t_1)=0$, we have, from (\ref{eq2})
\begin{equation}\label{eq3.11}
H_{\alpha,\beta} ^{w}(X; t_1,t_2)=\frac{1}{\beta-\alpha}\log\int_{t_1}^{t_2}\left(\frac{xf(x)}{ F(t_2)}\right)^{\alpha+\beta-1}dx.
\end{equation}
Now we define
\begin{equation*}
M(t_2)=\int_{t_1}^{t_2}\left(\frac{xf(x)}{ F(t_2)}\right)^{\alpha+\beta-1}dx,
\end{equation*}
which on differentiation with respect to $t_2$, gives
\begin{equation}\label{eq3.111}
 M^{'}(t_2)=\frac{f(t_2)}{F^{\alpha+\beta}(t_2)}\left[t_2^{\alpha+\beta-1}f^{\alpha+\beta-2}(t_2)F(t_2)-(\alpha+\beta-1)\int_{t_1}^{t_2}x^{\alpha+\beta-1}f^{\alpha+\beta-1}(x)dx\right].
\end{equation}
The above shows that the sign of $ M^{'}(t_2)$ relies only upon the factor in square braces in (\ref{eq3.111}). Again, define
\begin{equation*}
 N(t_2)=t_2^{\alpha+\beta-1}f^{\alpha+\beta-2}(t_2)F(t_2)-(\alpha+\beta-1)\int_{t_1}^{t_2}x^{\alpha+\beta-1}f^{\alpha+\beta-1}(x)dx.
\end{equation*}
Clearly $ N(t_2)|_{t_2=t_1}=0$ and it's derivative
\begin{equation*}
   N^{'}(t_2)= (\alpha+\beta-1)(t_2f(t_2))^{\alpha+\beta-2}F(t_2)-(\alpha+\beta-2)t_2^{\alpha+\beta-1}f^{\alpha+\beta-3}(t_2)(f^{2}(t_2)-f^{'}(t_2)F(t_2))\geq0,
\end{equation*}
if $\alpha+\beta<2$ and $F(x)$ is log-concave. Thus, $N(t_2)$ is positive, which in turn gives that $H_{\alpha,\beta}^{w}(X; t_1,t_2)$ is an increasing function of $t_2$. \\
(ii) Now we consider a random variable $X_{t_1}$ with distribution function $G(x)=\frac{F(x)-F(t_1)}{1-F(t_1)}$. Clearly, $G(t_1)=0$ and $G(x)$ is twice differentiable. When $F(x)$ is log-concave, $G(x)$ is easily verified to also be log-concave. Thus, applying Theorem \ref{th3.11}(i) to $X_{t_1}$, we find $H_{\alpha,\beta}^{w}(X_{t_1}; t_1,t_2)$ is increasing in $t_2$. At the same time, it is simple to verify that
\begin{equation*}
  H_{\alpha,\beta}^{w}(X_{t_1}; t_1,t_2)=H_{\alpha,\beta} ^{w}(X; t_1,t_2).
\end{equation*}
Thus, $H_{\alpha,\beta}^{w}(X; t_1,t_2)$ is also an increasing function of $t_2$. Due to symmetry, the log-concavity also implies that $H_{\alpha,\beta}^{w}(X; t_1,t_2)$ is decreasing in $t_1$. Hence the result follows.$\hfill\square$\\

\hspace*{.2in} We would like to remark that many commonly used distributions have log-concave cumulative distribution functions. For example, exponential, Pareto, lognormal, power distribution, Weibull distribution with shape parameter in $(0,1)$, gamma distribution with shape parameter in $(0,1)$ etc. are log-concave. This shows a wide range of applicability of the above result.\\
\hspace*{.2in} In the sequel we obtain some bounds for $H_{\alpha,\beta}^{w}(X; t_1,t_2)$. For the sake of brevity, the proofs are omitted.
\begin{t1}\label{th3.1}
For an absolutely continuous nonnegative random variable $X$, if $H_{\alpha,\beta} ^{w}(X; t_1,t_2)$ is increasing (decreasing) in $t_1$ for fixed $t_2$, then
\begin{equation*}
H_{\alpha,\beta}^{w}(X; t_1,t_2)\geq(\leq)\frac{1}{\beta-\alpha}\log\left[\frac{t_1^{\alpha+\beta-1}h_1^{\alpha+\beta-2}(t_1,t_2)}{(\alpha+\beta-1)}\right].
\end{equation*}
\end{t1}

\begin{t1}\label{th3.2}
If $H_{\alpha,\beta} ^{w}(X; t_1,t_2)$ is increasing (decreasing) in $t_2$ for fixed $t_1$, then
\begin{equation*}
H_{\alpha,\beta}^{w}(X; t_1,t_2)\leq(\geq)\frac{1}{\beta-\alpha}\log\left[\frac{t_2^{\alpha+\beta-1}h_2^{\alpha+\beta-2}(t_1,t_2)}{(\alpha+\beta-1)}\right].
\end{equation*}
\end{t1}
\hspace*{.2in} In the following theorem we give bound for $H_{\alpha,\beta}^{w}(X; t_1,t_2)$ based on monotonicity of $h_1(t_1,t_2)$ and $h_2(t_1,t_2)$.
\begin{t1}
Let $X$ be an absolutely continuous nonnegative random variable with density function $f(x)$ and distribution function $F(x)$. Then for $\alpha+\beta>(<)2$\\
(i) increasing $h_1(t_1,t_2)$ in $t_1$ implies
\begin{equation*}
H_{\alpha,\beta} ^{w}(X; t_1,t_2)\geq \left(\frac{\alpha+\beta-1}{\beta-\alpha}\right)\log\left[ t_1 h_1(t_1,t_2)\right];
\end{equation*}
(ii) decreasing $h_2(t_1,t_2)$ in $t_2$ implies
\begin{equation*}
H_{\alpha,\beta} ^{w}(X; t_1,t_2)\geq\left(\frac{\alpha+\beta-1}{\beta-\alpha}\right)\log\left[ t_2 h_2(t_1,t_2)\right].
\end{equation*}
\end{t1}
Proof: (i) By recalling (\ref{eq2}),
\begin{eqnarray*}
H_{\alpha,\beta} ^{w}(X;t_1,t_2)&=&\frac{1}{\beta-\alpha}\log\int_{t_1}^{t_2}\left(x\frac{f(x)}{F(t_2)-F(x)}\frac{F(t_2)-F(x)}{F(t_2)-F(t_1)}\right)^{\alpha+\beta-1}dx\\
&=&\frac{1}{\beta-\alpha}\log\int_{t_1}^{t_2} x^{\alpha+\beta-1} h_1^{\alpha+\beta-1}(x,t_2)\left(\frac{F(t_2)-F(x)}{ F(t_2)-F(t_1)}\right)^{\alpha+\beta-1}dx.\\
\end{eqnarray*}
Since $\frac{F(t_2)-F(x)}{ F(t_2)-F(t_1)}\geq0$ and $x>t_1$ implies that $h_1(x,t_2)\geq h_1(t_1,t_2)$ then, we have
\begin{eqnarray*}
H_{\alpha,\beta} ^{w}(X;t_1,t_2) &\geq&\frac{1}{\beta-\alpha}\log\int_{t_1}^{t_2} t_1^{\alpha+\beta-1}h_1^{\alpha+\beta-1}(t_1,t_2)\left(\frac{F(t_2)-F(x)}{ F(t_2)-F(t_1)}\right)^{\alpha+\beta-1}dx \\
 &=& \frac{1}{\beta-\alpha}\left[\log\left(t_1^{\alpha+\beta-1} h_1^{\alpha+\beta-1}(t_1,t_2)\right)+\log\int_{t_1}^{t_2}\left(\frac{F(t_2)-F(x)}{ F(t_2)-F(t_1)}\right)^{\alpha+\beta-1}dx\right] \\
 &\geq&\frac{\alpha+\beta-1}{\beta-\alpha}\log\left(t_1 h_1(t_1,t_2)\right).
\end{eqnarray*}
The proof of the second part is similar.$\hfill\square$\\

\hspace*{.2in} Consider the following example in support of the above bounds.
\begin{ex}\label{exc342i}
For an absolutely continuous random variable $X$ having $f(x)=\frac{1}{b-a},~a<x<b,~a,b>0$, the GFR functions are $h_i(t_1,t_2)=1/(t_2-t_1),~~i=1,2.$
It is not very difficult to see that $h_1(t_1,t_2)$ is increasing in $t_1$ and $h_2(t_1,t_2)$ is decreasing in $t_2$. Now, from Table 1, we have
\begin{equation*}
H_{\alpha,\beta}^{w}(X; t_1,t_2) =\frac{1}{\beta-\alpha}\log\left[\frac{1}{\alpha+\beta}\left(\frac{t_2^{\alpha+\beta}-t_1^{\alpha+\beta}}{(t_2-t_1)^{\alpha+\beta-1}}\right)\right]
\end{equation*}
\begin{equation*}
{\rm and,}~~~\frac{\alpha+\beta-1}{\beta-\alpha}\log[t_ih_i(t_1,t_2)]=\frac{\alpha+\beta-1}{\beta-\alpha}\log\left[\frac{t_i}{t_2-t_1}\right],~~i=1,2.
\end{equation*}
Then
$$d_i(t_1,t_1)=H_{\alpha,\beta} ^{w}(X; t_1,t_2)-\frac{\alpha+\beta-1}{\beta-\alpha}\log[t_ih_i(t_1,t_2)]\geq0,~~i=1,2$$
as shown in Figure \ref{fgc34.5iii}.
\begin{figure}[ht]
\centering
\begin{minipage}[b]{0.4\linewidth}
\includegraphics[height=4.5cm]{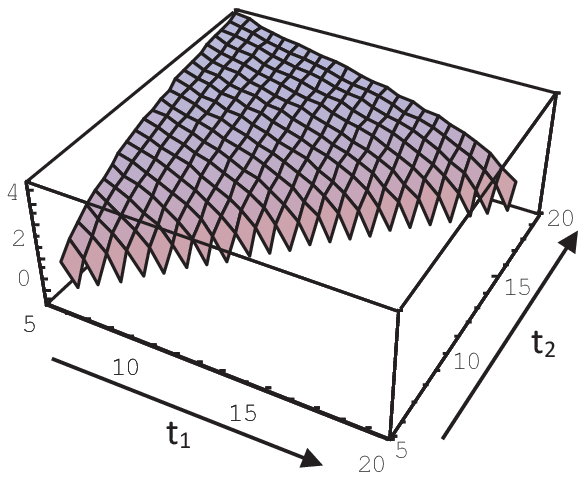}
\centering{$(i)$~Plot of $d_1(t_1,t_2)$ for $\alpha=0.5$ and $\beta=1.2$ against $(t_1, t_2)\in(5,20)$}
\end{minipage}
\quad
\begin{minipage}[b]{0.4\linewidth}
\includegraphics[height=4.5cm]{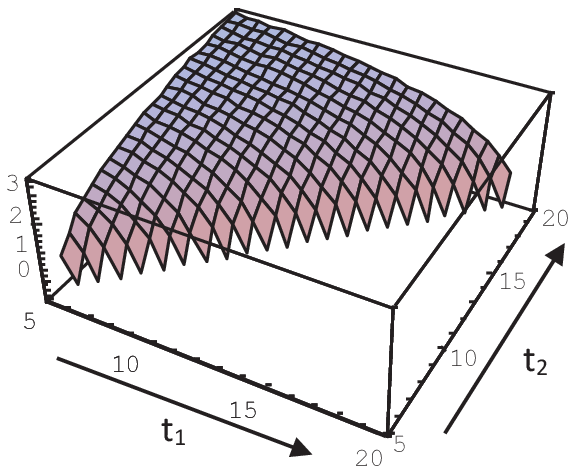}
\centering{$(ii)$~Plot of $d_2(t_1,t_2)$ for $\alpha=1.2$ and $\beta=2$ against $(t_1, t_2)\in(5,20)$}
\end{minipage}
\caption{\label{fgc34.5iii}Graphical representation of $d_1(t_1,t_2)$ and $d_2(t_1,t_2)$ (Example \ref{exc342i})}
\end{figure}\begin{figure}[ht]
\end{figure}
\end{ex}
\hspace*{.2in} Now we have a result which is applicable to large class of distributions that have monotone densities. Examples include exponential, Pareto and generalized Pareto, mixture of exponential, mixture of Paretos, Gamma and Weibull with shape parameters less than unity, folded symmetric distributions, to mention a few. The proof is omitted.
\begin{t1}\label{th4.5}
Let $X$ be an absolutely continuous nonnegative random variable. If $f(x)$ is increasing in $x>0$, then for all $t_1,t_2\in D$ and $\alpha+\beta>(<)2$\\
\begin{equation*}
H_{\alpha,\beta} ^{w}(X; t_1,t_2)\leq(\geq)\left(\frac{1}{\beta-\alpha}\right)\log E(X^{\alpha+\beta-1}|t_1<X<t_2)+\left(\frac{\alpha+\beta-2}{\beta-\alpha}\right)\log h_2(t_1,t_2),
\end{equation*}
\begin{equation*}
H_{\alpha,\beta} ^{w}(X; t_1,t_2)\geq(\leq)\left(\frac{1}{\beta-\alpha}\right)\log E(X^{\alpha+\beta-1}|t_1<X<t_2)+\left(\frac{\alpha+\beta-2}{\beta-\alpha}\right)\log h_1(t_1,t_2).
\end{equation*}
If $f(x)$ is decreasing in $x$, then the above inequalities are reversed.
\end{t1}

\hspace*{.2in} Now, we consider some inequalities based on WGIE.
\begin{p1}\label{pr3.1}
For an absolutely continuous nonnegative random variable $X$,
\begin{equation*}
H_{\alpha,\beta}^{w}(X; t_1,t_2)\leq\frac{1}{\alpha-\beta}\left(1-\int_{t_1}^{t_2}\left(\frac{xf(x)}{F(t_2)-F(t_1)}\right)^{\alpha+\beta-1}dx\right).
\end{equation*}
\end{p1}
\begin{t1}\label{th4.6}
Let $X$ be an absolutely continuous nonnegative random variable with distribution function $F(x)$ and $t_1,t_2\in D$. Then,
\begin{equation*}
H_{\alpha,\beta} ^{w}(X; t_1,t_2)\geq \frac{\alpha+\beta-1}{\beta-\alpha}\int_{t_1}^{t_2}\frac{f(x)}{F(t_2)-F(t_1)}\log x dx+\left(\frac{2-\alpha-\beta}{\beta-\alpha}\right)H(X; t_1,t_2)
\end{equation*}
where $H(X; t_1,t_2)=-\int_{t_1}^{t_2}\frac{f(x)}{ F(t_2)-F(t_1)}\log\frac{f(x)}{ F(t_2)-F(t_1)}dx$ is the interval Shannon entropy (cf. Misagh and Yari, 2011).
\end{t1}
Proof: From log-sum inequality, we have
\begin{eqnarray}\label{eq3.1}
\int_{t_1}^{t_2}f(x)\log\left(\frac{f(x)}{\left(\frac{xf(x)}{(F(t_2)-F(t_1)}\right)^{\alpha+\beta-1}}\right)dx&\geq&\left(\int_{t_1}^{t_2}f(x)dx\right)\log\left(\frac{\int_{t_1}^{t_2}f(x)dx}{\int_{t_1}^{t_2}\left(\frac{xf(x)}{F(t_2)-F(t_1)}\right)^{\alpha+\beta-1}dx}\right)\nonumber\\
&=&(F(t_2)-F(t_1))\left[\log(F(t_2)-F(t_1))\right.\nonumber\\
&&\left.-(\beta-\alpha)H_{\alpha,\beta} ^{w}(X; t_1,t_2)\right].
\end{eqnarray}
The left hand side of (\ref{eq3.1}) is
\begin{eqnarray}\label{eq3.2}
&&\int_{t_1}^{t_2}f(x)\log f(x)dx-(\alpha+\beta-1)\int_{t_1}^{t_2}f(x)\log\left(\frac{xf(x)}{F(t_2)-F(t_1)}\right)dx\nonumber\\
&&=(2-\alpha-\beta)\int_{t_1}^{t_2}f(x)\log f(x)dx\nonumber\\
&&-(\alpha+\beta-1)\left(\int_{t_1}^{t_2}f(x)\log x dx-(F(t_2)-F(t_1))\log(F(t_2)-F(t_1))\right)
\end{eqnarray}
From (\ref{eq2}), (\ref{eq3.1}) and (\ref{eq3.2}) we get the required result.$\hfill\square$\\

\hspace*{.2in} The following example illustrates the above theorem.
\begin{ex}\label{ex5}
Let $X$ be a random lifetime having $f(x)=2x$, $0<x<1$. Then
\begin{equation}\label{eq3.3}
H_{\alpha,\beta} ^{w}(X; t_1,t_2)=\frac{1}{\beta-\alpha}\log\left(\frac{2}{t_2^{2}-t_1^{2}}\right)^{\alpha+\beta-1}\left[\frac{t_2^{2(\alpha+\beta)-1}-t_1^{2(\alpha+\beta)-1}}{2(\alpha+\beta)-1}\right],
\end{equation}
\begin{equation}\label{eq3.4}
\frac{\alpha+\beta-1}{\beta-\alpha}\int_{t_1}^{t_2}\frac{f(x)}{F(t_2)-F(t_1)}\log x dx=\frac{\alpha+\beta-1}{(\beta-\alpha)(t_2^{2}-t_1^{2})}\left[t_2^{2}\log t_2-t_1^{2}\log t_1\right]-\frac{\alpha+\beta-1}{2(\beta-\alpha)},
\end{equation}
\begin{equation}\label{eq3.5}
H(X;t_1,t_2)=-\log2-\frac{1}{(t_2^{2}-t_1^{2})}\left[t_2^{2}\log t_2-t_1^{2}\log t_1\right]+\frac{1}{2}+\log(t_2^{2}-t_1^{2}).
\end{equation}
Now, from (\ref{eq3.3}), (\ref{eq3.4}) and (\ref{eq3.5}), we have for $\alpha+\beta>(<)2$
\begin{eqnarray*}
c(t_1,t_2)&\stackrel{{\rm def}}{=}&\frac{1}{\beta-\alpha}\log\left(\frac{2}{t_2^{2}-t_1^{2}}\right)^{\alpha+\beta-1}\left[\frac{t_2^{2(\alpha+\beta)-1}-t_1^{2(\alpha+\beta)-1}}{2(\alpha+\beta)-1}\right]\\
&&-\left(\frac{(\alpha+\beta-1)}{(\beta-\alpha)(t_2^{2}-t_1^{2})}\left[t_2^{2}\log t_2-t_1^{2}\log t_1\right]-\frac{\alpha+\beta-1}{2(\beta-\alpha)}\right)\\
&&-\frac{2-\alpha-\beta}{\beta-\alpha}\left(-\log2-\frac{1}{(t_2^{2}-t_1^{2})}\left[t_2^{2}\log t_2-t_1^{2}\log t_1\right]+\frac{1}{2}+\log(t_2^{2}-t_1^{2})\right)\geq0,
\end{eqnarray*}
as shown in Figure \ref{fig3.3}, satisfying Theorem \ref{th4.6}.
\begin{figure}[ht]
\centering
\begin{minipage}[b]{0.4\linewidth}
\includegraphics[height=4cm]{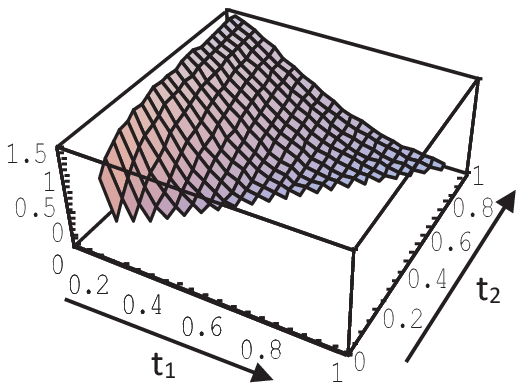}
\centering{$(i)$~Plot of $c(t_1,t_2)$ for $\alpha+\beta>2$ against $t_1,t_2\in(0,1)$}
\end{minipage}
\quad
\begin{minipage}[b]{0.4\linewidth}
\includegraphics[height=4cm]{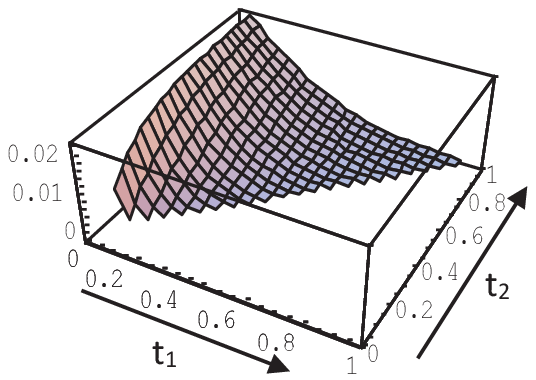}
\centering{$(ii)$ Plot of $c(t_1,t_2)$ for $\alpha+\beta<2$ against $t_1,t_2\in(0,1)$}
\end{minipage}
\caption{\label{fig3.3}Graphical representation of $c(t_1,t_2)$ for $\alpha+\beta>(<)2$ (Example \ref{ex5})}
\end{figure}
\end{ex}
\hspace*{.2in} We conclude this section by addressing an important question: does WGIE determine the distribution uniquely? Theorem \ref{th2.1} provides an answer to the same. Recalling (\ref{eq2}), we have
\begin{equation}\label{eq2.3}
e^{\left({\beta-\alpha}\right)H_{\alpha,\beta} ^{w}(X; t_1,t_2)}=\int_{t_1}^{t_2}\left(\frac{xf(x)}{ F(t_2)-F(t_1)}\right)^{\alpha+\beta-1}dx.
\end{equation}
 Differentiating (\ref{eq2.3}) with respect to $t_1$, we get
\begin{equation}\label{eq2.4}
(\beta-\alpha) \frac{\partial}{\partial t_1}H_{\alpha,\beta} ^{w}(X; t_1,t_2)= -(t_1h_1(t_1,t_2))^{\alpha+\beta-1}e^{-\left({\beta-\alpha}\right)H_{\alpha,\beta} ^{w}(X; t_1,t_2)}+(\alpha+\beta-1) h_1(t_1,t_2).
\end{equation}
Similarly, differentiating (\ref{eq2.3}) with respect to $t_2$, we get
\begin{equation}\label{eq2.5}
(\beta-\alpha) \frac{\partial}{\partial t_2}H_{\alpha,\beta} ^{w}(X; t_1,t_2)= (t_2h_2(t_1,t_2))^{\alpha+\beta-1}e^{-\left({\beta-\alpha}\right)H_{\alpha,\beta} ^{w}(X; t_1,t_2)}-(\alpha+\beta-1) h_2(t_1,t_2).
\end{equation}
Then, for any fixed $t_1$ and arbitrary $t_2$, $h_1(t_1,t_2)$ is a positive solution of the equation $\eta(x_{t_2})=0$, where
\begin{equation*}
\eta(x_{t_2})\stackrel{{\rm def}}{=} -(t_1x_{t_2})^{\alpha+\beta-1}e^{-\left({\beta-\alpha}\right)H_{\alpha,\beta} ^{w}(X; t_1,t_2)}+(\alpha+\beta-1) x_{t_2}-(\beta-\alpha) \frac{\partial}{\partial t_1}H_{\alpha,\beta} ^{w}(X; t_1,t_2).
\end{equation*}
Similarly, for any fixed $t_2$ and arbitrary $t_1$, $h_2(t_1,t_2)$ is a positive solution of the equation $\zeta(y_{t_1})=0$, where
\begin{equation*}
\zeta(y_{t_1})\stackrel{{\rm def}}{=} -(t_2y_{t_1})^{\alpha+\beta-1}e^{-\left({\beta-\alpha}\right)H_{\alpha,\beta} ^{w}(X; t_1,t_2)}+(\alpha+\beta-1) y_{t_1}+(\beta-\alpha) \frac{\partial}{\partial t_2}H_{\alpha,\beta} ^{w}(X; t_1,t_2).
\end{equation*}
Differentiating $\eta(x_{t_2})$ and $\zeta(y_{t_1})$ with respect to $x_{t_2}$ and $y_{t_1}$, respectively, we get
\begin{equation*}
\frac{\partial}{\partial x_{t_2}}\eta(x_{t_2})=(\alpha+\beta-1)-(\alpha+\beta-1) t_1^{\alpha+\beta-1}(x_{t_2})^{\alpha+\beta-2}e^{-\left({\beta-\alpha}\right)H_{\alpha,\beta} ^{w}(X; t_1,t_2)}
\end{equation*}
\begin{equation*}
{\rm and,}~~\frac{\partial}{\partial y_{t_1}}\zeta(y_{t_1}) =-(\alpha+\beta-1)t_2^{\alpha+\beta-1}(y_{t_1})^{\alpha+\beta-2}e^{-\left({\beta-\alpha}\right)H_{\alpha,\beta} ^{w}(X; t_1,t_2)}+(\alpha+\beta-1).
\end{equation*}
Furthermore, we consider second order derivative of $\eta(x_{t_2})$ and $\zeta(y_{t_1})$ with respect to $x_{t_2}$ and $y_{t_1}$, given by
\begin{equation*}
\frac{\partial^{2}}{\partial x_{t_2}^{2}}\eta(x_{t_2})=-(\alpha+\beta-1)(\alpha+\beta-2)t_1^{\alpha+\beta-1}(x_{t_2})^{\alpha+\beta -3}e^{-\left({\beta-\alpha}\right)H_{\alpha,\beta} ^{w}(X; t_1,t_2)}
\end{equation*}
\begin{equation*}
{\rm  and,}~~\frac{\partial^{2}}{\partial y_{t_1}^{2}}\zeta(y_{t_1})=-(\alpha+\beta-1)(\alpha+\beta-2)t_2^{\alpha+\beta-1}(y_{t_1})^{\alpha+\beta-3}e^{-\left({\beta-\alpha}\right)H_{\alpha,\beta} ^{w}(X; t_1,t_2)}.
\end{equation*}
Now, $\frac{\partial}{\partial x_{t_2}}\eta(x_{t_2})=0$ gives $x_{t_2}=\left[t_1^{-(\alpha+\beta-1)}e^{(\beta-\alpha)H_{\alpha,\beta} ^{w}(X; t_1,t_2)}\right]^{\frac{1}{\alpha+\beta-2}}=x_{t_2}^{o}$, say and $\frac{\partial}{\partial y_{t_1}}\zeta(y_{t_1})=0$ gives $y_{t_1}=\left[t_2^{-(\alpha+\beta-1)}e^{(\beta-\alpha)H_{\alpha,\beta}^{w}(X; t_1,t_2)}\right]^{\frac{1}{\alpha+\beta-2}}=y_{t_1}^{o}$, say.
\begin{t1}\label{th2.1}
For an absolutely continuous nonnegative random variable $X$, if $H_{\alpha,\beta}^{w}(X; t_1,t_2)$ is increasing in $t_1$ (for fixed $t_2$) and decreasing in $t_2$ (for fixed $t_1$), then\\
$(i)$ $\eta(x_{t_2})=0$ and $\zeta(y_{t_1})=0$ have unique solutions $x_{t_2}=h_1(t_1,t_2)$ and $y_{t_1}=h_2(t_1,t_2)$ if $\eta(x_{t_2}^{0})=0$ and $\zeta(y_{t_1}^{0})=0$. Thus, the distribution is determined uniquely;\\
$(ii)$ $\eta(x_{t_2})=0$ (resp. $\zeta(y_{t_1})=0$) has two solutions if $\eta(x_{t_2}^{0})\neq0$ (resp. $\zeta(y_{t_1}^{0})\neq0$). Of these two solutions, at least one should be $h_1(t_1,t_2)$ (resp. $h_2(t_1,t_2)$).
\end{t1}
Proof: We prove the theorem in two different cases.\\
{\bf Case 1:} Let $\alpha+\beta<2$ then $\eta(0)=-(\beta-\alpha) \frac{\partial}{\partial t_1}H_{\alpha,\beta}^{w}(X; t_1,t_2)<0$ and $\eta(\infty)=\infty$,  since $H_{\alpha,\beta} ^{w}(X; t_1,t_2)$ is increasing in $t_1$ (for fixed $t_2$) and decreasing in $t_2$ (for fixed $t_1$). Also, $\eta(x_{t_2})$ is a convex function with minimum occurring at $x_{t_2}=x_{t_2}^{0}$. Thus $\eta(x_{t_2})=0$ has a unique solution. Further, $\zeta(0)=(\beta-\alpha) \frac{\partial}{\partial t_2}H_{\alpha,\beta} ^{w}(X; t_1,t_2)<0$ and $\zeta(\infty)=\infty$. Again, $\zeta(y_{t_1})$ is a convex function with minimum occurring at $y_{t_1}=y_{t_1}^{0}$. So, $\zeta(y_{t_1})=0$ has a unique solution.\\
{\bf Case 2:} Let $\alpha+\beta>2$ then $\eta(0)=-(\beta-\alpha) \frac{\partial}{\partial t_1}H_{\alpha,\beta}^{w}(X; t_1,t_2)<0$ and $\eta(\infty)=-\infty$. Further, one can see that $\eta(x_{t_2})$ is a concave function with maximum occurring at $x_{t_2}=x_{t_2}^{0}$. Therefore $\eta(x_{t_2})=0$ has a unique solution if $\eta(x_{t_2}^{0})=0$. Also, $\zeta(0)=(\beta-\alpha) \frac{\partial}{\partial t_2}H_{\alpha,\beta}^{w}(X; t_1,t_2)<0$, $\zeta(\infty)=-\infty$ and $\zeta(y_{t_1})$ is a concave function with maximum occurring at $y_{t_1}=y_{t_1}^{0}$. Thus, $\zeta(y_{t_1})=0$ has a unique solution when $\zeta(y_{t_1}^{0})=0$.\\
\hspace*{.2in} Therefore, both the equations $\eta(x_{t_2})=0$ and $\zeta(y_{t_1})=0$ have unique positive solutions $h_1(t_1,t_2)$ and $h_2(t_1,t_2)$, respectively, if $\eta(x_{t_2}^{0})=0$ and  $\zeta(y_{t_1}^{0})=0$. Hence $H_{\alpha,\beta}^{w}(X; t_1,t_2)$ uniquely determines the GFR functions which in turn determines the distribution function uniquely (cf. Navarro and Ruiz, 1996).
\begin{r1}
It can be shown that an analogous result also holds if $H_{\alpha,\beta} ^{w}(X; t_1,t_2)$ is decreasing in $t_1$ (for fixed $t_2$) and increasing in $t_2$ (for fixed $t_1$).
\end{r1}
\begin{r1} Note that $\eta(x_{t_2})=0$ and  $\zeta(y_{t_1})=0$ have unique solution for all $(t_1,t_2)\in D$ when $\eta(x_{t_2}^{0})=0$ and  $\zeta(y_{t_1}^{0})=0$ which gives that $H_{\alpha,\beta}^{w}(X; t_1,t_2)=\frac{1}{\beta-\alpha}\log\frac{1}{\alpha+\beta}\left[\frac{t_2^{\alpha+\beta}-t_1^{\alpha+\beta}}{(t_2-t_1)^{\alpha+\beta-1}}\right]$, i.e., the WGIE of uniform distribution over $(a,b)$, $a<b$. Thus the uniform distribution over $(a,b)$, $a<b$ can be characterized by $H_{\alpha,\beta} ^{w}(X; t_1,t_2)=\frac{1}{\beta-\alpha}\log\frac{1}{\alpha+\beta}\left[\frac{t_2^{\alpha+\beta}-t_1^{\alpha+\beta}}{(t_2-t_1)^{\alpha+\beta-1}}\right]$.
\end{r1}
\hspace*{.2in} Below we give one example where both the solutions of $\eta(x_{t_2})=0$ are GFR as claimed in Theorem \ref{th2.1}.
\begin{ex}\label{ex6}
Let $X$ follow beta distribution with density function $f(x)=2x, 0\leq x\leq1$. Then for $\alpha=1.8$ and $\beta=2.5$, $H_{\alpha,\beta} ^{w}(X; t_1,t_2)=\frac{1}{0.7}\log\left[\left(\frac{2^{3.3}}{7.6}\right)\left(\frac{t_{2}^{7.6}-t_{1}^{7.6}}{(t_2^{2}-t_1^{2})^{3.3}}\right)\right]$ is increasing in $t_1$ (for fixed $t_2$) and decreasing in $t_2$ (for fixed $t_1$) for $t_1,t_2\in (0,1)$. Then, we have
\begin{equation}\label{55}
\frac{x_{t_2}^{0}}{h_1(t_1,t_2)}=t_1^{-4.3}\left(\frac{2}{t_2^{2}-t_1^{2}}\right)^{2.3}\left[\frac{t_{2}^{7.6}-t_{1}^{7.6}}{7.6}\right]^{\frac{1}{2.3}}=\nu(t_1,t_2),~say
\end{equation}
\begin{figure}
\centering
\includegraphics[width=6cm,keepaspectratio]{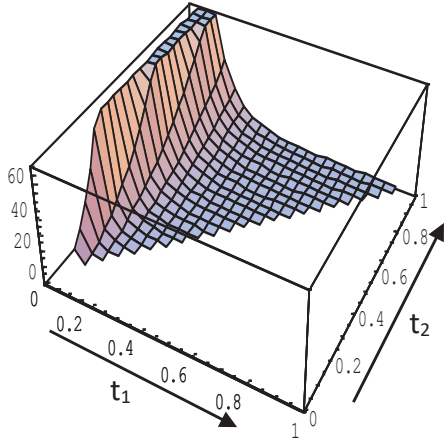}
\vspace{.5cm} \caption{\label{fig3.4} Plot of $\left[\nu(t_1,t_2)-1\right]$ against $t_1\in(0,1)$ and $t_2\in(0,1)$ (Example \ref{ex6})}
\end{figure}
which is greater than 1 for all $t_1,t_2\in (0,1)$ as shown in Figure \ref{fig3.4}. So, for every $t_1,t_2>0$, $\eta(x_{t_2})=0$ has two positive solutions as $h_1(t_1,t_2)$ and $h_1^{\ast}(t_1,t_2)$ such that $h_1(t_1,t_2)<x_{t_2}^{0}<h_1^{\ast}(t_1,t_2)$ and therefore $h_1^{\ast}(t_1,t_2)$ must be a GFR.
\end{ex}
\section{Simulation study and analysis of a real data set}
In this section we estimate WGIE and further carried out a simulation study to illustrate the performance of the estimator using simulated and real data sets. All the simulation works have been done using R-software.
\subsection{Simulation study}
Here we estimate WGIE by using Monte-Carlo simulation study and examine the performance of estimated values of $H_{\alpha,\beta}^{w}(X; t_1,t_2)$ i.e., $\widehat H_{\alpha,\beta}^{w}(X; t_1,t_2)$. To this aim, we use the method of maximum likelihood. Let $X$ follow Exp($\lambda$). First we estimate the unknown parameter $\lambda$ i.e., $\widehat\lambda$ by using maximum likelihood estimation method and then use it in (\ref{eq2}) to get the maximum likelihood estimator (MLE) of $H_{\alpha,\beta}^{w}(X; t_1,t_2)$ which is given by
$$\widehat H_{\alpha,\beta}^{w}(X; t_1,t_2)=\frac{1}{\beta-\alpha}\log\int_{t_1}^{t_2}\left(\frac{\widehat\lambda x\exp(-\widehat\lambda x)}{\exp(-\widehat\lambda t_1)-\exp(-\widehat\lambda t_2)}\right)^{\alpha+\beta-1}dx.$$
To illustrated the performance of the estimator, we generate samples from double truncated exponential distribution with parameter value 2. The estimated values are computed based on 1000 simulations each of size $n$ ($n=50,100,500,1000$) for different truncation limits and $\alpha+\beta<(>)2$. Averages are calculated from these 1000 values of $\widehat H_{\alpha,\beta}^{w}(X; t_1,t_2)$ which give their final values. Bias and mean squared error (MSE) of $\widehat H_{\alpha,\beta}^{w}(X; t_1,t_2)$ are also calculated. In Table 2-3, we present the estimates, bias and MSE for $\alpha+\beta<(>)2$, respectively. It is clear from Table 2 that $\widehat H_{\alpha,\beta}^{w}(X; t_1,t_2)$ increases as $t_1$ and $t_2$ increases (when the other is fixed) for $\alpha+\beta<2$. Also, from Table 3 we observe that $\widehat H_{\alpha,\beta}^{w}(X; t_1,t_2)$ increases with respect to $t_1$ and decreases with respect to $t_2$ (when the other is fixed) for $\alpha+\beta>2$. It is worthwhile to remark that this outcome is in accordance with the monotonicity of $H_{\alpha,\beta}^{w}(X; t_1,t_2)$ for Exp(2). The results of simulation studies show that as the sample size increases, absolute values of bias and MSE decreases and for large sample estimates are almost unbiased.
\subsection{Analysis of real data set}
In this subsection, we further analyze a real data set. Here we consider the data set representing the times of successive failures of the air conditioning system of each member of a fleet of Boeing 720 jet airplanes which was analyzed by Proschan (1963). For illustrative purpose, we consider a single airplane namely, Plane 7912. The hours of flying time between successive failures for this plane are given below.\\
Data Set (Plane 7912): 1, 3, 5, 7, 11, 11, 11, 12, 14, 14, 14, 16, 16, 20,
21, 23, 42, 47, 52, 62, 71, 71, 87, 90, 95, 120, 120, 225, 246, 261.\\
As has been observed by Proschan (1963), the exponential distribution with hazard rate $\lambda$ can be fitted to this data set. We verify the same through a goodness-of-fit test. The Kolmogorov-Smirnov (K-S) distance between the empirical distribution and the fitted distribution functions and the associated $p$-value were obtained as 0.1581 and 0.5602, respectively. Now we obtain the estimates of $H_{\alpha,\beta}^{w}(X; t_1,t_2)$. To this aim, we estimate $\lambda$ on using the method of maximum likelihood for different truncation limits and then used them to find $\widehat H_{\alpha,\beta}^{w}(X; t_1,t_2)$. Table 4 provides the estimated values of WGIE for different truncation limit $(t_1,t_2)$ and $\alpha+\beta<(>)2$. It is clear from Table 4, that $\widehat H_{\alpha,\beta}^{w}(X; t_1,t_2)$ increases with respect to $t_1$ for $\alpha+\beta<(>)2$ and also increases with respect to $t_2$ for $\alpha+\beta<2$  but, interestingly, decreases with respect to $t_2$ for $\alpha+\beta>2$ (when the other is fixed). Therefore, the monotonic behavior of the estimates as observed for simulated data are validated by the airplane data set as well.

\section{Application of WGIE in reliability modeling}
\begin{figure}[ht]
\centering
\begin{minipage}[b]{0.4\linewidth}
\includegraphics[height=4.6cm]{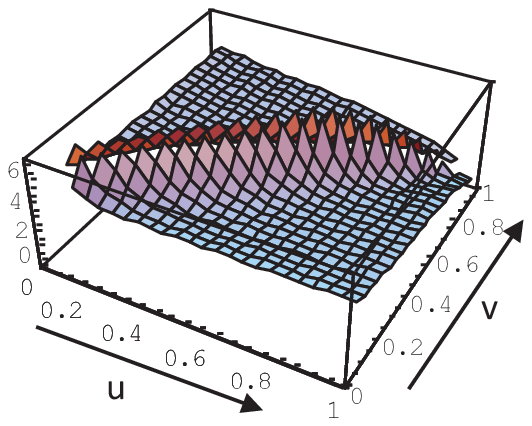}
\centering{$(i)$ Plot of $\kappa_{\alpha,\beta} ^{w}( u,v)$ for $\alpha=1.5$ and $\beta=2$}
\end{minipage}
\quad
\begin{minipage}[b]{0.4\linewidth}
\includegraphics[height=4.6cm]{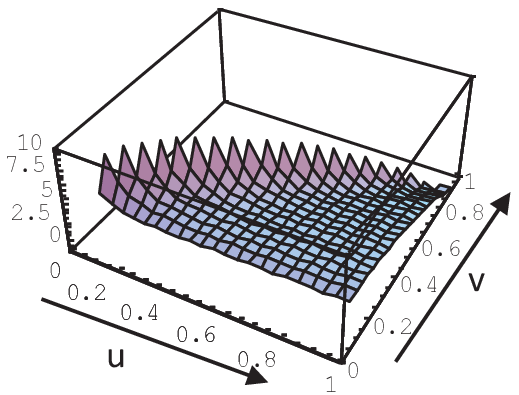}
\centering{$(ii)$ Plot of $\eta_{\alpha,\beta} ^{w}(u,v)$ for $\alpha=1.5$ and $\beta=2$}
\end{minipage}
\caption{\label{fig5.1}Graphical representation of $\kappa_{\alpha,\beta} ^{w}( u,v)$ and $\eta_{\alpha,\beta} ^{w}(u,v)$ for $\alpha+\beta>2$}
\end{figure}

\begin{figure}[ht]
\centering
\begin{minipage}[b]{0.4\linewidth}
\includegraphics[height=4.6cm]{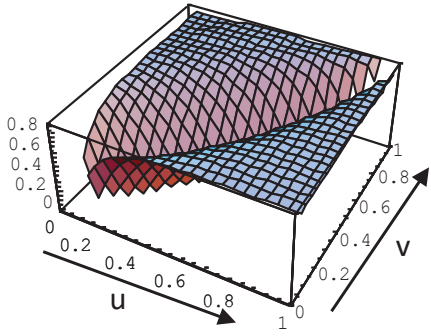}
\centering{$(i)$ Difference of WGIE between EE and Gamma distribution for $\alpha=1.5$ and $\beta=2$}
\end{minipage}
\quad
\begin{minipage}[b]{0.4\linewidth}
\includegraphics[height=4.6cm]{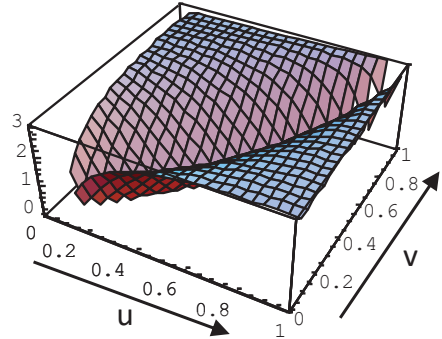}
\centering{$(ii)$ Difference of WGIE between EE and Weibull distribution  for $\alpha=1.5$ and $\beta=2$}
\end{minipage}
\caption{\label{fig5.2}Graphical representation for difference of $H_{\alpha,\beta}^{w}(\cdot; t_1,t_2)$}
\end{figure}

In this section, we investigate the role of our proposed measure (WGIE) in reliability modeling for a real-life problem. We know that entropy is the measure of uncertainty (randomness) of a process or system and the probability distribution which best represents the current state of knowledge for the given data set is one with maximum entropy. The principle of Maximum Entropy enunciated by Jaynes (1957) is a technique that can be used to estimate input probability more generally. It states that out of all distributions consistent with a given set of constraints choose one that maximizes entropy. For some flavour of fascinating growth of maximum entropy model and information theoretic approach for model selection, one may refer to Kapur (1994) and Burnham and Anderson (2003), respectively. According to the information-theoretic approach for model selection due to Burnham and Anderson (2003), for a given data set, the best fitted model is the one which has maximum entropy associated with it. Between two models, the more accurate model will be the one with larger entropy. To this aim various extensions of Shannon entropy have been proposed in the literature that may have more information (uncertainty) about a distribution than the information given by the Shannon entropy. In order to take into account the qualitative characteristic of information, the WGIE can be used for comparing different probabilistic models when we do not know the actual probability distribution that generated some data.\\
\hspace*{.2in} To see the effectiveness of $H_{\alpha,\beta}^{w}(X; t_1,t_2)$ in reliability modeling we consider the data set arose in tests on endurance of deep groove ball bearings (Lawless, 1986, P. 228). The observations are the number of million revolutions before failure for each 23 ball bearings in the life test; the individual bearings were inspected periodically to determine whether failure had occurred. The data set are given bellow.\\
{\bf Data Set:}  17.88, 28.92, 33.00, 41.52, 42.12, 45.60, 48.80, 51.84, 51.96, 54.12, 55.56, 67.80, 68.64, 68.64, 68.88, 84.12, 93.12, 98.64, 105.12, 105.84, 127.92, 128.04 and 173.40.\\
Gupta and Kundu (2001) fitted the following three distributions to analyze the data set.\\
$(i)$ The Gamma distribution
\begin{equation}\label{eq5.1}
f(x)=\frac{\lambda^{a}}{\Gamma a}(x)^{a-1}e^{-\lambda x}; \quad a,~\lambda,~x>0
\end{equation}
with $\widehat a=4.0196$ and $\widehat \lambda=0.0556$.\\
$(ii)$ The Weibull distribution
\begin{equation}\label{eq5.2}
f(x)=a\lambda(\lambda x)^{a-1}e^{-(\lambda x)^{a}}; \quad a,~\lambda,~x>0
\end{equation}
with $\widehat a=2.1050$ and $\widehat \lambda=0.0122$.\\
$(iii)$ The exponentiated exponential (EE) distribution
\begin{equation}\label{eq5.3}
f(x)=a\lambda(1-e^{-\lambda x})^{a-1}e^{-\lambda x}; \quad a,~\lambda,~x>0
\end{equation}
with $\widehat a=5.2589$ and $\widehat \lambda=0.0314$.\\
They have claimed that for the given data set EE distribution provides a better fit compared to Weibull or Gamma distributions.\\
\hspace*{.2in} We now examine the role of WGIE for comparing statistical models to be fitted to the given data set. Let $X$ be a nonnegative random variable which follow EE, Weibull and Gamma distribution as given in (\ref{eq5.3}), (\ref{eq5.2}) and (\ref{eq5.1}), respectively. Then Figure \ref{fig5.1} shows that $H_{\alpha,\beta} ^{w}(X; t_1,t_2)-H(X; t_1,t_2)=\kappa_{\alpha,\beta} ^{w}(X; t_1,t_2)$, say and $H_{\alpha,\beta} ^{w}(X; t_1,t_2)-H^{w}(X; t_1,t_2)=\eta_{\alpha,\beta} ^{w}(X; t_1,t_2)$, say are positive for some $(\alpha,\beta)$ so that $\alpha+\beta>2$ when $X$ follows EE distribution enabling $H_{\alpha,\beta}^{w}(X; t_1,t_2)$ has larger uncertainty, in view of qualitative characteristic of information, than $H(X; t_1,t_2)$ and $H^{w}(X; t_1,t_2)$.
We also plot difference of two weighted generalized interval entropies in which first one follows EE distribution and other follows Gamma or Weibull distribution. It is shown in Figure \ref{fig5.2} that the differences are always positive for the same values of $\alpha$ and $\beta$. Note that the substitutions $t_1=-\log u$ and $t_2=-\log v$ have been used while plotting curves so that $\kappa_{\alpha,\beta}^{w}(X; t_1,t_2)=\kappa_{\alpha,\beta} ^{w}(u,v)$, say.
Though various entropy measures are available in the literature, one should choose that entropy measure which has maximum uncertainty associated with a distribution. In agreement with Gupta and Kundu (2001) if $X$ follows (\ref{eq5.3}), i.e., the best fitted model for the given data set, then one can see that the uncertainty contained in $H_{\alpha,\beta}^{w}(X; t_1,t_2)$ is more than $H(X; t_1,t_2)$ and $H^{w}(X; t_1,t_2)$ for some specific values of parameters and uncertainty contained in $H_{\alpha,\beta}^{w}(\cdot; t_1,t_2)$ for EE distribution is greater than for Gamma or Weibull distributions.\\
\hspace*{.2in} Hence from the above discussion, we can conclude that when the qualitative characteristic of information is taken into consideration, the WGIE contains more (average) information than the interval Shannon entropy and weighted interval entropy, respectively. The uncertainty content in $H_{\alpha,\beta}^{w}(X; t_1,t_2)$ for EE distribution is more as compare to Gamma and Weibull distributions which indeed enable one to find out the best fitted model. More work is needed in this direction.

\section*{Acknowledgements}
The financial support (Ref. No. 2/48(4)/2015/NBHM(R.P.)/R$\&$D II/14130\&11629) rendered by the NBHM, Department of Atomic Energy, Government of India is gratefully acknowledged.

\begin{landscape}
\begin{table}[h]\label{tbl2.1}
{
\caption{WGIE for some well-known distributions.}
\begin{center}
 \begin{tabular}{|c|c|c|}

   \hline
   Model & $f(x)$ & $H_{\alpha,\beta}^{w}(X; t_1,t_2)$ \\
    \hline

 Uniform &$\frac{1}{b-a};~ a<x<b$ &$\frac{1}{\beta-\alpha}\log\left[\frac{1}{\alpha+\beta}\left(\frac{t_2^{\alpha+\beta}-t_1^{\alpha+\beta}}{(t_2-t_1)^{\alpha+\beta-1}}\right)\right]$ \\
  Exponential &$\theta\exp(-\theta x);~ x,\theta>0$ &$\frac{1}{\beta-\alpha}\log\left[\frac{1}{\theta(\alpha+\beta-1)^{\alpha+\beta}}\left(\frac{\gamma (\alpha+\beta,\theta(\alpha+\beta-1)t_2)
-\gamma(\alpha+\beta,\theta(\alpha+\beta-1)t_1)}{\left(e^{-\theta t_1}-e^{-\theta t_2}\right)^{\alpha+\beta-1}}\right)\right]$ \\
 Power&$\frac{b}{a}(\frac{x}{a})^{b-1};~ 0<x<a,~b>0$ &$\frac{1}{\beta-\alpha}\log\left[\frac{b^{\alpha+\beta-1}}{a^{b(\alpha+\beta-1)}(b(\alpha+\beta-1)+1)}\left(\frac{t_2^{b(\alpha+\beta-1)+1}-t_1^{b(\alpha+\beta-1)+1}}{((\frac{t_2}{a})^{b}-(\frac{t_1}{a})^{b})^{\alpha+\beta-1}}\right)\right]$ \\
 Beta&$cx^{c-1};~ 0<x<1,~c>0$ &$\frac{1}{\beta-\alpha}\log\left[\frac{c^{\alpha+\beta-1}}{((c-1)(\alpha+\beta-1)+1)}\left(\frac{t_2^{(c-1)(\alpha+\beta-1)+1}-t_1^{(c-1)(\alpha+\beta-1)+1}}{(t_2^{c}-t_1^{c})^{\alpha+\beta-1}}\right)\right]$ \\
 Pareto I &$\frac{b}{a}\left(\frac{x}{a}\right)^{-(b+1)};~x>a$ &$\frac{1}{\beta-\alpha}\log\left[\frac{(ba^{b})^{\alpha+\beta-1}}{1-b(\alpha+\beta-1)}\left(\frac{t_2^{1-b(\alpha+\beta-1)}-t_1^{1-b(\alpha+\beta-1)}}{\left((\frac{a}{t_1})^{b}-(\frac{a}{t_2})^{b}\right)^{\alpha+\beta-1}}\right)\right]$ \\
 Gamma&$\frac{b^{n}x^{n-1}e^{-nx}}{\Gamma n}$;~$b,n>0$ &$\frac{1}{\beta-\alpha}\log\left[\frac{b^{n(\alpha+\beta-1)}}{(n(\alpha+\beta-1))^{n(\alpha+\beta-1)+1}}\left(\frac{\gamma (n(\alpha+\beta-1)+1,n(\alpha+\beta-1)t_2)}{(\gamma (n,b t_2)
-\gamma(n,b t_1))^{\alpha+\beta-1}}-\frac{\gamma(n(\alpha+\beta-1)+1,n(\alpha+\beta-1)t_1)}{(\gamma (n,b t_2)
-\gamma(n,b t_1))^{\alpha+\beta-1}}\right)\right]$ \\
   \hline

 \end{tabular}
 \end{center}
  }
 \end{table}

\begin{table}[h]\label{tbl4.1}
{\tiny
\caption{$\widehat H_{\alpha,\beta}^{w}(X;t_1,t_2)$, Bias and MSE for $\alpha=0.5$ and $\beta=1.2$
($n=50,100,500,1000$).}
\begin{center}
\begin{tabular}{|c|c|c|c|}\hline $(t_1,t_2)$ &$\widehat H_{\alpha,\beta}^{w}(X;t_1,t_2)$
&Bias &MSE\\
\hline
\begin{tabular}{c}
   \\\hline
  (1,3)\\
(1,5)\\
(1,7)\\
(3,11)\\
(5,11)\\
(7,11)\\
\end{tabular}
& \begin{tabular}{cccc}
                        $n=50$ & $n=100$ & $n=500$ & $n=1000$ \\
                        \hline
                        0.5819471&0.5819471&0.5819471&0.5819471\\
                        0.7106222&0.7198491& 0.7162191&0.7182855\\
                         0.7259325&0.7239873&0.7281197&0.728812\\
                        1.516652&1.517459&1.520471&1.519483\\
                        1.948473&1.953715&1.953276&1.952496\\
                        2.242633&2.246337&2.246826&2.248049
                        \end{tabular}
& \begin{tabular}{cccc}
                        $n=50$ & $n=100$ & $n=500$ & $n=1000$ \\
                        \hline
                        -0.006758232&-0.008428941&-0.001759892&-0.000549343\\
                      -0.007044602&0.002182274&-0.001447732&0.0006187335\\
                         -0.002964676&-0.004909903&-0.000777495&-8.52E-05\\
                        -0.003035095&-0.002227601&0.000784352&-0.00020357\\
                        -0.004291961&0.000949689&0.000510775&-0.000269345\\
                        -0.004948353&-0.00124425&-0.000755268&0.000467877
                       \end{tabular}
& \begin{tabular}{cccc}
                        $n=50$ & $n=100$ & $n=500$ & $n=1000$ \\
                        \hline
                        0.005720018&0.002837238&0.000557682&0.000270947\\
                       0.01206394&0.005801127&0.001100345&0.0005762391\\
                        0.01326288&0.006749558&0.001407349&	0.000643114\\
                        0.006999292&0.003566514&0.000797305&0.000374912\\
                        0.005702713&0.002912027&0.000627575&0.000297053\\
                        0.004983525&0.00219102&0.000459015&0.000239054
                       \end{tabular}\\
\hline
\end{tabular}
\end{center}
}
\end{table}

\begin{table}[h]\label{tbl4.2}
{\tiny
\caption{$\widehat H_{\alpha,\beta}^{w}(X;t_1,t_2)$, Bias and MSE for $\alpha=1.5$ and $\beta=2$
($n=50,100,500,1000$).}
\begin{center}
\begin{tabular}{|c|c|c|c|}\hline $(t_1,t_2)$ &$\widehat H_{\alpha,\beta} ^{w}(X;t_1,t_2)$
&Bias &MSE\\
\hline
\begin{tabular}{c}
   \\\hline
  (1,3)\\
(1,5)\\
(1,7)\\
(3,11)\\
(5,11)\\
(7,11)\\
 \end{tabular}
& \begin{tabular}{cccc}
                        $n=50$ & $n=100$ & $n=500$ & $n=1000$ \\
                        \hline
                       1.399301&1.365262&1.359219&1.355347\\
                       1.313549&1.283728&1.26998&1.270621\\
                       1.288068&1.26814&1.266727&1.264182\\
                       6.102378&6.089941&6.083496&6.076423\\
                       8.523865&8.505244&8.491835&8.497547\\
                       10.1491&10.13756&10.11984&10.12334
                      \end{tabular}
& \begin{tabular}{cccc}
                        $n=50$ & $n=100$ & $n=500$ & $n=1000$ \\
                        \hline
                        0.04279845&	0.008759404&0.002716423&-0.001155106\\
                       0.04678905&0.0169678&0.003219549&0.003860197\\
                        0.02295513&0.003027008&0.001613672&-0.00093077\\
                        0.02487513&0.0124379&0.005992474&-0.001080494\\
                        0.02807229&0.009451012&-0.003957554&0.001754433\\
                        0.02723631&0.01569853&-0.002016842&0.001483488
                        \end{tabular}
& \begin{tabular}{cccc}
                        $n=50$ & $n=100$ & $n=500$ & $n=1000$ \\
                        \hline
                        0.08189&0.03469952&0.007786411&0.003519915\\
                        0.0892493&0.04152996&0.007746569&0.004271309\\
                         0.08354724&0.03897862&0.00780873&0.004158212\\
                        0.1546349&0.07190441&0.01471098&0.007160713\\
                        0.1553508&0.08132905&0.01705033&0.007974807\\
                        0.168429&0.08449432&0.01599597&0.008746537
                       \end{tabular}\\
\hline
\end{tabular}
\end{center}
}
\end{table}

\begin{table}[h]\label{tbl4.3}
{
\caption{Estimated values of $H_{\alpha,\beta}^{w}(X;t_1,t_2)$ for the air plane data for different truncation limits $(t_1,t_2)$ and $\alpha+\beta<(>)2$.}
\begin{center}
\begin{tabular}{|c|c|}\hline&$\widehat H_{\alpha,\beta}^{w}(X;t_1,t_2)$\\
\hline
\begin{tabular}{c}$(\alpha,\beta)\backslash(t_1,t_2)$
   \\\hline
  (0.5,1.2)\\
(1.5,2)\\
\end{tabular}

&\begin{tabular}{cccccc}
                        $(10,20)$ & $(10,50)$ & $(10,90)$ & $(12,200)$&$(50,200)$&$(90,200)$ \\
                        \hline
                       3.613115&4.487881&5.434379&6.153233&6.524566&6.572326\\
                       6.384205& 4.832763&4.803095&4.930697&8.432027&11.87169

                        \end{tabular}\\
\hline
\end{tabular}
\end{center}
}
\end{table}
\end{landscape}


\begin{thebibliography}{}
{\small
\bibitem{bd} Baig, M.A.K. and Dar, J.G. (2008), Generalized residual entropy function and its applications. {\it European Journal of Pure and Applied Mathematics}, {\bf1}, 30-40.
\bibitem{bg} Belis, M. and Guia\c{s}u, S. (1968), A quantitative-qualitative measure of information in cybernetic systems. {\it IEEE Transcations on Information Theory}, {\bf 14}, 593-594.
\bibitem{ba} Burnham, K.P. and Anderson, D. (2003), {\it Model Selection and Multimodel Inference: A Practical Information-Theoretic Approach}. 2nd Ed., Springer-Verlag, New York.
\bibitem{d} Das, S. (2017), On weighted generalized entropy. {\it Communications in Statistics- Theory \& Methods}, {\bf46(12)}, 5707-5727.
\bibitem{dl} Di Crescenzo, A. and Longobardi, M. (2006), On weighted residual and past entropies. {\it Scientiae Mathematicae
Japonicae}, {\bf64}, 255-266.
\bibitem{f} Fisher, R.A. (1934), The effect of methods of ascertainment upon the estimation of frequencies. {\it Annals
of Eugenics}, {\bf6}, 13-25.
\bibitem{gk} Gupta, R.D. and Kundu, D. (2001), Exponentiated exponential family: an alternative to gamma and Weibull distributions. {\it Biometrical Journal}, {\bf43(1)}, 117-130.
\bibitem{j} Jaynes, E.T. (1957), Information theory and statistical mechanics. {\it Physical Review}, {\bf106(4)}, 620-630.
\bibitem{j} Kapur, J.N. (1994), {\it Measures of Information and Their Applications}. Wiley-Interscience.
\bibitem{k2b} Kayal, S. (2015), On generalized dynamic survival and failure entropies of order $(\alpha,\beta)$. {\it Statistics and Probability Letters}, {\bf96}, 123-132.
\bibitem{kt} Kumar, V. and Taneja, H.C. (2011), Some characterization results on generalized cumulative residual entropy measure. {\it Statistics and Probability Letters}, {\bf81}, 1072-1077.
\bibitem{k} Kundu, C. (2015), Generalized measures of information for truncated random variables. {\it Metrika}, {\bf78(4)}, 415-435.
\bibitem{ks} Kundu, C. and Singh, S. (2020), On generalized interval entropy. {\it Communications in Statistics- Theory \& Methods}, {\bf49(8)}, 1989-2007.
\bibitem{l} Lawless, J.F. (1986), {\it Statistical Models and Methods for Lifetime Data}. Wiley, New York.
\bibitem{m} Minimol, S. (2017), On generalized dynamic cumulative past entropy measure. {\it Communications in Statistics- Theory \& Methods}, {\bf46(6)}, 2816-2822.
\bibitem{my} Misagh, F. and Yari, G.H. (2011), On weighted interval entropy. {\it Statistics and Probability Letters}, {\bf81}, 188-194.
\bibitem{nr} Navarro, J. and Ruiz, J.M. (1996), Failure rate functions for doubly truncated random variables. {\it IEEE Transactions on Reliability}, {\bf45(4)}, 685-690.
\bibitem{ny} Nourbakhsh, M. and Yari, G. (2017), Weighted renyi's entropy for lifetime distributions. {\it Communications in Statistics- Theory \& Methods}, {\bf46(14)}, 7085-7098.
\bibitem{p} Proschan, F. (1963), Theoretical explanation of observed decreasing failure rate. {\it Technometrics}, {\bf5(3)}, 375-383.
\bibitem{ranr} Rajesh, G., Abdul-Sathar, E.I., Reshmi, K.V. and Nair, K.R.M. (2014), Bivariate generalized cumulative residual entropy. {\it Sankhy$\bar{a}$: A}, {\bf76(1)}, 101-122.
\bibitem{ranr} Rajesh, G., Abdul-Sathar, E.I. and Rohini, S. Nair. (2017), On dynamic weighted survival entropy of order $\alpha$. {\it Communications in Statistics- Theory \& Methods}, {\bf46(5)}, 2139-2150.
\bibitem{r} Rao, C.R. (1965), {\it Linear Statistical Inference and its Applications}. Wiley, New York.
\bibitem{r} R$\acute{e}$nyi, A. (1961), On measures of entropy and information. {\it In: Proceeding of the Fourth Berkeley Symposium on Mathematical Statistics and Probability}, {\bf 1}, 547-561.
\bibitem{sbr} Sekeh, Y.S., Borzadaran, G.R.M. and Raknabadi, A.H.R. (2014), Some results based on a version of the generalized dynamic entropies. {\it Communications in Statistics- Theory \& Methods}, {\bf43(14)}, 2989-3006.
\bibitem{sc} Shangari, D. and Chen, J. (2012), Partial monotonicity of entropy measures. {\it Statistics and Probability Letters}, {\bf82(11)}, 1935-1940.
\bibitem{s} Shannon, C.E. (1948), A mathematical theory of communication. {\it Bell System Technical Journal}, {\bf 27}, 379-423, 623-656.
\bibitem{sk} Singh, S. and Kundu, C. (2019), On weighted Renyi's entropy for double truncated distribution. {\it Communications in Statistics- Theory \& Methods}, {\bf48(10)}, 2562-2579.
\bibitem{v} Varma, R.S. (1966), Generalizations of Renyi's entropy of order $\alpha$. {\it Journal of Mathematical Sciences}, {\bf1}, 34-48.
\bibitem{w} Wallis, G. (1996), Using spatio-temporal correlations to learn invariant object recognition. {\it Neural Networks}, {\bf9 (9)}, 1513-1519.
}
\end{thebibliography}
\end{document}